\documentclass[11pt]{article}
\usepackage{graphicx}
\usepackage{balance} 
\usepackage{lipsum} 
\usepackage[sc]{mathpazo} 
\usepackage[T1]{fontenc} 
\usepackage{microtype} 
\usepackage[hmarginratio=1:1,top=30mm,columnsep=15pt]{geometry} 
\usepackage{multicol} 
\usepackage[hang, small,labelfont=bf,up,textfont=it,up]{caption} 
\usepackage{booktabs} 
\usepackage{float} 
\usepackage{hyperref} 
\usepackage{lettrine} 
\usepackage{paralist} 
\usepackage{bbm}
\usepackage{float}
\usepackage{amsmath,amssymb}  
\usepackage{color} 
\usepackage[square]{natbib}
\usepackage{abstract} 

\usepackage{titlesec} 
\renewcommand\thesection{\Roman{section}} 
\renewcommand\thesubsection{\Roman{subsection}} 
\titleformat{\section}[block]{\large\scshape\centering}{\thesection.}{1em}{} 
\titleformat{\subsection}[block]{\large}{\thesubsection.}{1em}{} 

\usepackage{fancyhdr} 
\newtheorem{theo}{Theorem}

\newtheorem{coro}{Corollary}

\newtheorem{defi}{Definition}
\newtheorem{exa}{Example}

\newtheorem{prop}{Proposition}

\newtheorem{rem}{Remark}

\newcommand{\RR}{\mathbb{R}}

\newcommand{\PP}{\mathbb{P}}

\newcommand{\ZZ}{\mathbb{Z}}
\newcommand{\NN}{\mathbb{N}}

\newcommand{\EE}{\mathbb{E}}

\newcommand{\mc}{\mathcal}

\newcommand{\mb}{\mathbb}
\newcommand{\lfled}{\longrightarrow}
\newcommand{\fled}{\rightarrow}
\newcommand{\­}{\neq}

\newcommand{\²}{\leqslant}
\newcommand{\³}{\geqslant}
\newcommand{\Cov}{\mbox{Cov}}
\newcommand{\Lip}{\mbox{Lip}}
\newcommand{\Var}{\mbox{Var}}

\newcommand{\e}{\mathbf{e}}
\newcommand{\h}{\widehat}
\newcommand{\1}{\mb{1}}

\title{\vspace{-15mm}\fontsize{25pt}{10pt}\selectfont\textbf{Dependent Lindeberg CLT - Finite Dimensional for Empirical Processes of Cluster Functionals}\thanks{This research has been conducted as part of the project Labex MME-DII (ANR11-LBX-0023-01)}} 
\author{\large \textsc{Jos\'e-Gregorio G\'omez}\footnote{D\'epartement de Math\'ematiques, Universit\'e de Cergy - Pontoise. 95000 Cergy - Pontoise, France. E-mail address: \href{mailto:jose.gomez@u-cergy.fr}{jose.gomez@u-cergy.fr}}}
\vspace{-5mm}
\date{}

\begin{document}
\maketitle 
\thispagestyle{fancy} 


\begin{abstract}
Drees and Rootz\'en [2010] have proven central limit theorems (CLT) for empirical processes
of extreme values cluster functionals built from $\beta$-mixing processes. The problem with this family of $\beta$-mixing processes is that it is quite restrictive, as has been shown by Andrews [1984]. We expand this result to a more general dependent processes family, known as weakly dependent processes in the sense of Doukhan and Louhichi [1999], but in finite-dimensional convergence (fidis). We show an example where the application of the CLT-fidis is sufficient in several cases, including a small simulation of the extremogram introduced by Davis and Mikosch [2009] to confirm the efficacy of our result.
\vspace{0.2 cm}
\begin{quotation}
\end{quotation}
\end{abstract}



\section{Introduction}

In light of recent developments in massive data processing via {\it parallel processing}, it is convenient to consider the construction of statistics in function of data blocks. In the case of extremes (rare events), we have very little data that is relevant to our estimations, but instead they are hidden among a large mass of "common data". Thus comes the natural idea of considering clustering of extremes, which here consists of obtaining the smaller sub-block of extreme values on each block, while conveniently suppressing "common" data in each block , generally assigned a null value. Such null values may be mathematically inoffensive, yet computationally they are an obstacle when it is our aim to obtain quick results. These and many other reasons encourage the study of extremes cluster functionals. \ This paper aims to offer a small contribution to the asymptotic behaviour of extremes cluster functionals. More particularly, an extension of the dependence condition of \citet{Drees2010}'s central limit theorem in finite-dimensional convergence (CLT-fidis) for the empirical processes of extremes cluster functionals.  

In order to do this, we use mainly \citet{Bardet2007}'s dependent Lindeberg method which addresses the construction of CLTs for dependent processes under the usual Lindeberg condition, if the sequence $T = T(n)$ (which summarizes the dependence of the process) tends towards zero when the number of random variables $n$ is large. Particularly, this term of $T$-dependence can be written as a sum of covariances which can be bounded by weak dependence coefficients defined by \citet{Doukhan1999}. Therefore, for weakly dependent processes with convenient decrease rates in the weak-dependence coefficients, we obtain CLT-fidis for empirical processes of cluster functionals (EPCFs). 

Several reasons motivate this extension. The main one is that weak dependence is a very general property including certain non-mixing processes: e.g. take a AR(1)-input, solution of the recursion 
\begin{align}\label{AR}
X_k=\dfrac{1}{b}\big(X_{k-1}+\xi_k\big),\quad k\in\ZZ,
\end{align}
 where $b\³2$ is an integer and $(\xi_k)_{k\in\NN}$ are independent and uniformly distributed random variables on the set $U(b):=\{0,1,\ldots, b-1\}$. This process is not mixing in the sense of Rosenblatt, as this is shown in \citep{Andrews}  for $b=2$ and  in \citep{Ango} for $b>2$, however \citet{Doukhan1999} proved that such a process is weakly dependent. More generally, under weak conditions, all the usual causal or non causal time series are weakly dependent processes: this is the case for instance of Gaussian, associated, linear, ARCH($\infty$), etc. 

This document is organised as follows: Section II offers the definition of EPCFs, generalized by \citet{Drees2010} for the multidimensional case and developed first by \citet{Yun2000}  and \citet{Segers2003} for the real case. In Section III we provide a general CLT-fidis for these empirical processes through \citet{Bardet2007}'s Lindeberg method, we define weak dependence and we provide some examples of weakly dependent processes. Finally, we apply the initial theorem to this type of dependent processes. In Section IV we develop an example where the finite - dimensional convergence of the EPCFs is sufficient; this is, the estimator of the \citet{Davis2009}'s extremogram. Proofs are given in Section V.

\section{Empirical Processes of Extremes Cluster Functionals}

In this section we outline some necessary basic definitions and hypotheses that we will consider throughout this document in order to prove limit theorems of empirical processes of extreme cluster functionals. Roughly, an extreme cluster functional is a map that works on blocks (arbitrary-length but not random-length) of "extreme" random variables in such a way that the map remains invariant under extreme clusters, which in this paper is the smaller sub-block that contains all the extreme values (and the non-extreme values among them) of the given block. Besides, this application is null when there are no extremes within the block. 

We will mention some examples in a cursory way just to better understand the definitions. Some technical details of these examples will be shown in Section IV. 
 
 \subsection{Cluster Functionals}
Let $(E,\mc{E})$ be a measurable subspace of $(\RR^d,\mc{B}(\RR^d))$ for some $d\³1$ such that $0\in E$. Following the deterministic definition of  \citet{Drees2010}\footnote{This definition is given by  \cite{Yun2000} and \cite{Segers2003}, for the real case}, we consider the set of $E$-valued sequences of finite length, {\it i.e.}, 
$$E_\cup:=\{(x_1, \ldots, x_r): x_i\in E  \quad\forall i=1,\ldots, r; \quad\forall r\in\NN \},$$
equipped with the $\sigma$-field $\mc{E}_\cup$ induced by Borel-$\sigma$-fields  on $E^r$, for $r\in\NN$. Then, if  $x\in E_\cup$, we can write $x=(x_1, \ldots, x_r)$ for some $r\in\NN$. The \textbf{core}\footnote{ Note that the core also considers the null values that exist between the non-null values. \\Ex. $(0,1,2,0,0,3,0,1,0,0)^c=(1,2,0,0,3,0,1)$, which is the smaller sub-block of $x=(0,1,2,0,0,3,0,1,0,0)$ which contains all non-null values as well as the null values between them.} $x^c\in E_\cup$ of  $x$ is defined by 
 \[x^c:=\left\{\begin{array}{cc}
 (x_{r_I}, x_{r_I+1}, \ldots, x_{r_S}), & \mbox{ if  $x\­ 0_r$ (the null element in $E^r$)}    \\ \\
 0,  & \mbox{  otherwise}
\end{array}\right.\]
where $r_I: =\min\{i\in \{1,\ldots,r\}: x_i \­ 0\}$ (first non-null value of the block $x$) and $r_{S}: =\max\{i\in \{1,\ldots,r\}: x_i \­ 0\}$ (last non-null value of the block $x$). A \textbf{cluster functional} is a measurable map $f: (E_\cup, \mc{E}_\cup)\lfled (\RR, \mc{B}(\RR))$ such that 
 \begin{align}\label{fc}
 f(x)=f(x^c),\qquad\mbox{ for all }x\in E_\cup,\ \ \mbox{ and } \ f(0_r)=0\  (\forall r\³1).
 \end{align}
  
Under the properties (\ref{fc}), it is easy to build a large amount of examples of cluster functionals. Nevertheless, the typical examples used to build estimators through these cluster functionals are functionals of the type: 
\begin{align}\label{combi}
 f(x_1, \ldots, x_r)=\sum_{i=1}^r \phi(x_i),
\end{align}
where $\phi: E\lfled \RR$ is such that $\phi(0)=0$. Generally speaking, these functions $\phi$  are  indicator functions (or functions which are product of another measurable function $H:E \lfled \RR$  with an indicator function). Another classic example is the component-wise maximum of a cluster:
 \begin{align}\label{combi2}
 f(x_1, \ldots, x_r)=\max_{1\²i\²r} x_i, 
\end{align}
for $E=[0,\infty)$.

\subsection{Empirical Processes of Cluster Functionals}

Now, we want to apply cluster functionals to blocks of $E$-valued random variables excesses over a determined thresholds sequence and to define the empirical process indexed by these functionals. 

Let us consider $E$-valued normalized random variables $(X_{n,i})_{1\²i\²n, n\in \NN}$, defined on some probability space $(\Omega, \mc{A}, \PP)$, which are row-wise stationary, this is, $(X_{n,i})_{1\²i\²n}$ is stationary for each $n\in \NN$. Here, those normalized random variables $X_{n,i}$ are built from another random process $(X_i)_{i\in\ZZ}$, in a way such that the normalization maps all non-extreme values to zero. Additionally, it should satisfy that the sequence of conditional distributions of $X_{n,1}$ given that $X_{n,1}$ belongs of the failure set $A \subseteq E\setminus \{0\}$ ({\it i.e. }$P_n(\cdot| A):=\PP\{X_{n,1}\in \cdot  |X_{n,1}\in A\}$),  converge weakly to some non-degenerate limit. 

For instance,  for a real-valued stationary $(X_i)_{i\in\NN}$ with marginal cumulative distribution function $F$, let $(u_n)_{n\in\NN}$ be a non-decreasing sequence of thresholds such that $u_n \uparrow x_F$, where $$x_F=\sup\{x\in\RR: F(x)<1\},\qquad v_n=\PP\{X_1>u_n\} \underset{n\to\infty}{\lfled} 0.$$ Quote that the tail distribution function of $X_i$ may be  asymptotically degenerated , which means that there exists a point $a\in\RR$ such that $$\overline{P}_n(x)=\PP\{X_1-u_n>x|X_1>u_n\}\underset{n\to\infty}{\lfled} \1_{x\²a}.$$  
However, if $F$ belongs to the  domain of attraction of some extreme-value distribution, then by a result in  \citep{Pickands1975}, there exists $\gamma\in\RR$ and a sequence of positive constants $(a_n)_{n\in\NN}$ (depending on the sequence $u_n$) such that
\[P_n(x)=\PP\{X_{n,1}>x|X_1>u_n\} \underset{n\to\infty}{\lfled} \left\{\begin{array}{ll}
(1+\gamma x)^{-1/\gamma}_+, & \mbox{if  $\gamma\­0$}     \\
 \e^{-x},  & \mbox{if $\gamma=0$}
\end{array}\right.\]
locally uniform in $(0,\infty)$, where 
\begin{align}\label{N1}
X_{n,i}=\left(\frac{X_i-u_n}{a_n}\right)_+:= \max\left\{ \frac{X_i-u_n}{a_n} , 0 \right\}, \qquad \text{ for }1\leqslant i \leqslant n.
\end{align}

For the multidimensional case, let $\mb{X}=(X_i)_{i\in\NN}$ be a stationary $\RR^d$-valued time series such that all components of $X_i$ have the same marginal distribution. Since such time series $\mb{X}$ may exhibit dependence across coordinates and over time, if $\| \cdot \|$ denotes an arbitrary norm on $\RR^d$, then an interesting normalization for study the extreme values of $\mb{X}$ would be:  
\begin{align}\label{norextremogram}
X_{n,i}=u_n^{-1} X_i \1 \{ \|X_i\| > u_n\},  \qquad \text{ for }1\leqslant i \leqslant n.
\end{align}
where $(u_n)_{n\in\NN}$ is a sequence of high quantiles of the process. 

\begin{defi}[Empirical Process of Cluster Functionals \citep{Drees2010}]\

Let $Y_{n,j}$ be the $j$-th block of $r_n$ consecutive values of the $n$-th row of $(X_{n,i})$. Thus there are $m_n:=[n/r_n]=\max\{j\in \NN: j \² n/r_n\} $ blocks $Y_{n,j}:=(X_{n,i})_{(j-1)r_n+1\²i\² jr_n}$ of length $r_n$,   with $1\² j \² m_n$.  Moreover, since $(X_{n,i})_{1\²i\²n}$ is stationary for each $n$, then we can denote by $Y_n$  to the ``generic block" such that $Y_n \overset{\mc{D}}{=}Y_{n,1}$.

Now, let $\mc{F}$ be a class of  cluster functionals. The \textbf{"empirical process $Z_n$ of cluster functionals''} in $\mc{F}$, is the process $(Z_n(f))_{f\in\mc{F}}$ defined by
\begin{align}\label{PE}
Z_n(f):= \frac{1}{\sqrt{n v_n}}\sum_{j=1}^{m_n}(f(Y_{n,j})-\EE f(Y_{n,j})),
\end{align}
where  $v_n:=\PP\{X_{n,1}\­ 0\}$.
\end{defi}

In order to begin approaching the convergence in fidis of the EPCF (\ref{PE}), observe that if the blocks $(Y_{n,j})_{1\²j\²m_n, n\in\NN}$ are independents and if we take in account the following essential convergence assumptions:
 \begin{itemize}
\item[\textbf{(C.1)}] $\EE\left[\left(f(Y_n)-\EE f(Y_n)\right)^2\mb{1}\left\{|f(Y_n)-\EE f(Y_n)|>\epsilon\sqrt{n v_n}\right\}\right]=o(r_n v_n)$, \\for all $\epsilon>0$, and  for all $f\in\mc{F}$.
\item[\textbf{(C.2)}] $(r_n v_n)^{-1}\Cov\left(f(Y_n),g(Y_n)\right)\lfled c(f,g)$, for all $f,g\in \mc{F}$,
\end{itemize}
with  $r_n \ll v_n^{-1} \ll n$, then the fidis of the empirical process  $(Z_n(f))_{f\in\mc{F}}$ of cluster functionals converge to the fidis of a Gaussian process $(Z(f))_{f\in\mc{F}}$ with the covariance function $c$.  

\citet{Drees2010} have proved CLTs for this process. In particular, they have proved a CLT-fidis of $(Z_n(f))_{f\in\mc{F}}$ using the  Bernstein blocks technique together with a  $\beta$-mixing coupling condition to boil down convergence to convergence of sums over i.i.d. blocks through \cite{Eberlein}'s technique involving the metric of total variation.

However, the family of mixing processes is quite restrictive. We can see this through a particularly simple example: the AR(1) - process defined in (\ref{AR}), which is not even $\alpha$-mixing. Therefore, the results in \citep{Drees2010} can not be used here. 
 
In our case, we will solve this problem in fidis through the Lindeberg method, developed by \citet{Bardet2007}, followed by its applications for weakly dependent random processes defined by \citet{Doukhan1999}. 

The weak spot under these weak dependence conditions, in the sense of \citet{Doukhan1999}, is that we have no coupling arguments to arrive to a uniform CLT, as \citet{Drees2010} have done in their paper by using the rich coupling properties of the $\beta$-mixing processes together with \citet{VdV}'s tightness criteria and asymptotic equicontinuity conditions.

The benefit of this work is that the convergence in fidis is sufficient in several examples and applications. Here, we will show a particular application in the Section IV. 

\section{Lindeberg Method and applications to CLT-fidis for empirical processes of cluster functionals}

In order to adapt \citet{Bardet2007}'s dependent Lindeberg method to CLT-fidis for empirical processes of cluster functionals, let us denote 
\begin{align}\label{N2}
W_{n,j}:=&W_{n,j}(f_1,\ldots, f_k)
\nonumber\\=&(n v_n)^{-1/2}\left(f_1(Y_{n,j})-\mb{E}f_1(Y_{n,j}), \ldots, f_k(Y_{n,j})-\mb{E}f_k(Y_{n,j})\right),
\end{align}
for $1\²j\²m_n$ and $(f_1, \ldots, f_k)\in\mc{F}^k$. Therefore, with this notation we have derived the following result:
\begin{theo}[Lindeberg CLT for cluster functionals]
Suppose that assumptions (C.1) and (C.2) hold with $r_n \ll v_n^{-1} \ll n$. Then, if  
\begin{align}\label{TDep}
T_t(m_n|f_1,\ldots,f_k):=\sum_{j=1}^{m_n}\left|\Cov(\mathbf{e}^{i<t,\sum_{s=1}^{j-1}W_{n,s}(f_1,\dots,f_k)>},\mathbf{e}^{i<t,W_{n,j}(f_1,\ldots,f_k)>})\right|\underset{n\fled\infty}{\lfled} 0
\end{align}
for all $t\in\RR^k$ and for all $k$-tuple of cluster functionals $(f_1,\ldots,f_k)\in\mc{F}^k$, the fidis of the empirical process  $(Z_n(f))_{f\in\mc{F}}$ of cluster functionals converge to the fidis of a Gaussian process $(Z(f))_{f\in\mc{F}}$ with the covariance function $c$.
\end{theo}

We have just seen that the convergence in fidis of $(Z_n(f))_{f\in\mc{F}}$ to a Gaussian law is obtained because $T_t(m_n|f_1,\ldots,f_k)$ converges to 0, for all $t\in\RR^k$ and for all $(f_1,\ldots, f_k)\in\mc{F}^k$ with $k\in\NN$. Actually this expression is related to the dependence of the random variables $(X_{n,i})_{1\²i\²n, n\in\NN}$. Note that $T_t(m_n|f_1, \ldots, f_k)$ is written in terms of sums of covariances, therefore using weak-dependence theory (see \citep{Dedecker2007}), we can give bounds for such $T_t(m_n|f_1,\ldots, f_k)$.  

\subsection{Weak Dependence} 
Let  $f: E^r\subseteq (\RR^d)^r \lfled \RR$ be a function, with $r\in\NN$.  As usual, we denote by:
$$\Lip (f):=\sup_{(x_1,\ldots,x_r)\­(y_1,\ldots,y_r)\in E^r}\dfrac{|f(x_1,\ldots,x_r)-f(y_1,\ldots,y_r)|}{\|x_1-y_1\| +\cdots + \|x_r-y_r\|}.$$
Similar to the definition of \citet{Doukhan1999}, we say that a triangular array of row-wise stationary $E$-valued random variables $\mathbb{M}=(X_{n,i})_{1\²i\²n, n\in\NN}$ is \textbf{$(\epsilon,\psi)$-weakly dependent} ($(\epsilon,\psi)$-WD) if there exist a function $\psi:(\NN)^2\times(\RR^+)^2\lfled \RR^+$, an infinite sequence of positive integers $(l_n)_{n\in\NN}$ with $l_n  \ll  n$, and a positive sequence  $(\epsilon_n(l_n))_{n\in\NN}$ decreasing to zero, such that 
\begin{align}\label{COV} 
\left|Cov\left(f(X_{n,i_1},\ldots,X_{n,i_u}),g(X_{n,j_1},\ldots,X_{n,j_v})\right)\right|\²\psi(u,v,Lip(f),Lip(g))\cdot\epsilon_n(l_n)
\end{align}
for all $(u,v)\in \NN\times \NN$, all $(i_1,\ldots,i_u)\in\NN^u$, $(j_1,\ldots,j_v)\in\NN^v$ with $i_1 <\cdots <i_u<i_u+l_n\²j_1<\cdots <j_v\²n$, and for all pair of functions  $(f,g)\in \Lambda^u(E)\times\Lambda^v(E)$, where $\Lambda^s(E):=\{ h:E^s\lfled \RR \text{ Lipschitzian with } \|h\|_\infty\²1 \text{ and }  \Lip(h)<\infty\}$.

\begin{rem}\rm Let us remember that $\mathbb{M}=(X_{n,i})_{1\²i\²n, n\in\NN}$ is constructed from another random process $\mathbb{X}=(X_i)_{i\in\ZZ}$. Therefore,  $\mathbb{M}$'s dependence properties are inherited from $\mathbb{X}$'s dependence properties. Even more so, if $\mathbb{X}$ is $(\epsilon,\psi)$-weakly dependent (in the usual sense defined by \citet{Doukhan1999} for random processes), then $\mathbb{M}$ is $(\tilde{\epsilon},\psi)$ - weakly dependent with $\tilde{\epsilon}_n(\cdot)=L_n\cdot \epsilon(\cdot)$, for some positive constant $L_n$ (which is written in function of  $\mathbb{M}$'s normalization constants). In this sense, if we want to study $\mathbb{M}$'s dependence properties, suffice it to take into account $\mathbb{X}$'s dependence properties. 
\end{rem}

We will consider four different particular cases of functions $\psi$ of weakly dependent processes:
\begin{enumerate}
\item If $\mathbb{X}$ is a causal random process, {\it i.e.} if there exist a function $H:D^\ZZ \lfled E $ and a $D$-valued sequence of independent and identically distributed random variables (i.i.d.r.v's) $(\xi_i)_{i\in\ZZ}$ such that $X_i=H(\xi_i,\xi_{i-1},\xi_{i-2},\ldots)$, for $i\³0$, is defined almost surely. 
Thus, the \textbf{$\theta$-weakly dependent causal condition} is defined by
\begin{align} \label{dp1}
\psi(u,v,Lip(f),Lip(g))=vLip(g).
\end{align}
In this case, we will simply denote $\theta(l)$ instead of $\epsilon(l)$.

\item If $\mathbb{X}$ is a non causal random process, the \textbf{$\eta$, $\kappa$, $\lambda$-weakly dependent conditions} are defined respectively by
\begin{align} 
		\psi(u,v,Lip(f),Lip(g))&=uLip(f)+vLip(g),\label{dp2}\\
		\psi(u,v,Lip(f),Lip(g))&=uvLip(f)Lip(g),\label{dp3}\\
		\psi(u,v,Lip(f),Lip(g))&=uLip(f)+vLip(g)+uvLip(f)Lip(g)\label{dp4}.
\end{align}
where we write $\eta(l)$, $\kappa(l)$ and $\lambda(l)$, respectively, instead of $\epsilon(l)$.
\end{enumerate}

\begin{exa}[Examples of Weakly Dependent Processes]\

\rm
Now, we give a little list of examples of weakly dependent processes with their dependence properties. Here, we consider $(\xi_n)_{n\in\ZZ}$ as a sequence of i.i.d.r.v's.
\begin{enumerate}
\item Suppose that $\mathbb{X}=(X_i)_{i\in\ZZ}$ is a ARMA($p,q$) - process, or more generally, a causal (respectively non causal) linear process such that $X_i=\sum_{j\³0}a_j\xi_{i-j}$ (respectively $X_i=\sum_{-\infty<j<\infty}a_j \xi_{i-j}$) for $i\in\ZZ$, where $a_j=\mc{O}(|j|^{-\nu})$ with  $\nu>1/2$. 

Then $\mathbb{X}$ is $\theta$- (resp. $\lambda$- ) weakly dependent with $\theta(l)=\lambda(l)=\mc{O}(l^{1/2-\nu})$. For more details, see \citep{Doukhan2002}.

In particular, the AR(1) - process (\ref{AR}) can be rewritten as the causal linear process $X_i=\sum_{j\³0}b^{-j-1} \xi_{i-j}$, with $\xi_0$ uniformly distributed on $\{0,\ldots, b-1\}$. In this case $X_0$ is uniformly distributed over $[0,1]$ and $\theta(l)\²b^{-l}$ (see Example 1  in \citep{Dedecker2004}).

\item Let $\mathbb{X}=(X_i)_{i\in\ZZ}$ be a GARCH($p,q$) - process or, more generally, a ARCH($\infty$) - process such that  $X_i=\sigma_i \xi_i$, where $\sigma_i^2=a_0+\sum_{j=0}^\infty a_j X_{i-j}^2$ for $i\in\ZZ$.
\begin{itemize}
\item For the GARCH($p,q$) case, if there is a constant $C>0$ and $\nu\in(0,1)$ such that for all $j\in\NN$, $0\²a_j\²C\nu^j$, then $\mathbb{X}$ is a $\lambda$-weakly dependent process such that $\lambda(l)=\mc{O}(\e^{-c\sqrt{l}})$ with $c>0$ (see \citep{Bardet2007}).
\item For the ARCH($\infty$) case, if there is a constant $C>0$ and $\nu>1$ such that for all $j\in\NN$, $0\²a_j\²Cj^{-\nu}$, then $\mathbb{X}$ is a $\lambda$-weakly dependent process with $\lambda(l)=\mc{O}(l^{-\nu+1})$ (see \citep{Doukhan2006}).
\end{itemize} 

\item Suppose that $\mathbb{X}=(X_i)_{i\in\ZZ}$ is a associated stationary process, then $\mathbb{X}$ is a $\lambda$-weakly dependent process such that $\lambda(l)=\mc{O}(\sup_{i\³l}\Cov(X_0,X_i))$.

\item If $\mathbb{X}=(X_i)_{i\in\ZZ}$ is a Gaussian process such that $\lim_{n\lfled \infty}\Cov(X_0,X_n)=0$. Then $\mathbb{X}$ is a $\lambda$-weakly dependent process with $\lambda(l)=\mc{O}(\sup_{i\³l}|\Cov(X_0,X_i)|)$. For details of the last two examples, see \citep{Doukhan1999}.
\end{enumerate}

Under suitable assumptions, the families of causal and non-causal bilinear processes, non-causal finite order Volterra processes, causal and non-causal infinite memory processes, etc., are also weakly dependent. For more details of weak-dependence properties of this processes, see the book of \cite{Dedecker2007}. 
\end{exa}

\subsection{Application of Theorem 1}

In the result below we give a CLT-fidis for cluster functionals of weakly dependent processes. The proof (Section V), we need that these functionals $f\in\mc{F}$ can be approximated through lipschitzian cluster functionals $(f_n)_ {n\³0}$. To built these $f_n$, we must consider certain truncation assumptions (C.4) on the functionals $f$ valued in the sub-blocks of length $r_n-l_n$ and a concentration condition (\ref{Concent}) on the probability measure. 

First, consider the following notation: if $Y=(x_1, x_2, \ldots, x_r)$, then  
\[Y^{(k)}=\left\{\begin{array}{cc}
 (x_1,\ldots, x_k) & \mbox{ if  $k \² r$}   \\ \\
 Y & \mbox{ if $k > r$.}
\end{array}\right.\]
Moreover, if $f\in\mc{F}$ is a cluster functional, then we denote $\Delta_n(f):=f(Y_n)-f(Y_n^{(r_n-l_n)})$, where $r_n$ is the length of the block $Y_n$ such that $l_n\ll r_n$.

\begin{prop} Suppose that (C.1), (C.2) and the following convergence conditions
\begin{itemize}
\item[\textbf{(C.3)}] $\sqrt{\Var \Delta_n(f)}=o\left(\dfrac{r^2_n}{n^2}\sqrt{n v_n}\right)$
\item[\textbf{(C.4)}] $\mb{E}^{1/2}|f(Y_n^{(r_n-l_n)})|^2 \mb{1}\left\{|f(Y_n^{(r_n-l_n)})|>\sqrt{n v_n}\right\}=o\left(\dfrac{r^2_n}{n^2}\sqrt{n v_n}\right)$,
\end{itemize}
are satisfied, with $r_n$, $l_n \underset{n\fled\infty}{\lfled}\infty $ such that $l_n \ll  r_n \ll v_n^{-1} \ll n$ and $r_n \ll \ln (n)$. Additionally, assume that the r.v's $(X_{n,i})_{1\²i\²n}$ are such that there exists positive real constants $C, \alpha, \rho$ ($\rho$ dependent of $n$) such that
\begin{align}\label{Concent}
\sup_{x\in E}\sup_{1\²i\²n}\PP\left\{ X_{n,i} \in B(x, \rho/2)\right\}\²C\rho^\alpha.
\end{align}
Then the fidis of the cluster functionals empirical process $(Z_n(f))_{f\in \mc{F}}$ converge to the fidis of a Gaussian process $(Z(f))_{f\in\mc{F}}$ with covariance function $c$ (defined in (C.2)), if the r.v's $(X_{n,i})_{1\²i\²n, n\in\NN}$ satisfies one of the following weak dependence cases:
\begin{enumerate}
\item[(D.1)] $\theta$-weakly dependent such that $\theta_n(k)=\mc{O}\left(k^{-\theta}\right)$ for some $\theta>0$ and $l_n^{-\theta}=o\left(r_n^2/n^2\right)$
\item[(D.2)] $\epsilon$-weakly dependent such that $\epsilon_n(k)=\mc{O}\left(k^{-\xi}\right)$ for some $\xi>0$ and $l_n^{-\xi}=o\left(r_n^3/n^3\right)$, where $\epsilon_n(\cdot)$ is any non causal weak dependence coefficient (\ref{dp2})-(\ref{dp4}).
\end{enumerate}
\end{prop}

Generally, the (C.2) convergence can be easily verified. However, through the following proposition (which is a similar result to \citet{Segers2003}'s Theorems 1 and 3) we provide some conditions, sufficient to verify (C.2) and which in some situations are easier to prove. Besides, this way we can give an alternative expression to the covariance  function $c$ (defined in (C.2)), as it is shown below in Corollary 1.
 
In order to carry this out, it is necessary to consider the following assumption:

\begin{itemize}
\item[\textbf{(TC)}] There is a sequence $W=(W_i)_{i\³1}$ of $E$-valued random variables such that, for all $k\in\NN$, the joint conditional distribution $$P_{(X_{n,i}, \mb{1}\{X_{n,i}=0\})_{1\²i\²k}| X_{n,1}\­0}$$ converges weakly to $P_{(W_i, \mb{1}\{W_i=0\})}$, and for all $f\in\mc{F}$ are a.s. continuous with respect to the distribution of $W^{(k)}=(W_1, \ldots, W_k)$ and $W^{(2:k)}=(W_2, \ldots, W_k)$ for all $k$, that is, 
$$\PP\{W^{(2:k)}\in D_{f, k-1}, W_i=0,\ \forall i>k \}=\PP\{W^{(k)}\in D_{f, k}, W_i=0,\ \forall i>k \}=0$$
where we denote by $D_{f,k}$  the set of discontinuities of $f|_{E^k}$.
\end{itemize}
\begin{rem}\rm
The existence of such sequence  $W$ is guaranteed in particular from Theorem 2 in \citep{Segers2003} with  $E=\RR$ and the normalization (\ref{N1}). There, Segers has shown that if $$\PP_{((X_{n,i})_{1\²i\²k}| X_1>u_n)}\underset{n\to\infty}{\lfled} -\log G_k,$$  where $G_k$ is some $k-$dimensional extreme value distribution for all $k\in\NN$, then there exists such "{\it tail chain}" $W=(W_i)_{i\in\NN}$ such that 
\begin{align}\label{tailchains}
P_{((X_{n,i}, \mb{1}\{X_{n,i}=0\})_{1\²i\²k}| X_1>u_n)} \overset{w}{\underset{_{n\to \infty}}{\lfled}} P_{(W_i, \mb{1}\{W_i=0\})_{1\²i\²k}},
\end{align} 
for all $k\in\NN$. 
\end{rem}
\begin{prop} Suppose that the r.v's $(X_{n,i})_{1\²i\²n, n\in\NN}$ satisfies one of the following weak dependence conditions:
\begin{enumerate}
\item[(D.1')] $\theta$-weakly dependent such that $\theta_n(l_n)=o\left(v_n^{p+1}\right)$
\item[(D.2')] $\eta$-weakly dependent such that $\eta_n(l_n)=o\left(v_n^{p+1}/r_n\right)$ 
\item[(D.3')] $\kappa$ ( resp. $\lambda$) -weakly dependent such that $\kappa_n(l_n)$ ( resp. $\lambda_n(l_n)$) $= o\left( v_n^{2p+1}/r_n\right)$,
\end{enumerate}
for some $p>0$, where $(r_n)$, $(l_n)$ are integer sequences such that $l_n \ll  r_n \ll v_n^{-1} \ll n$  with $l_n \underset{n\fled\infty}{\lfled}\infty$. Then, 
$$\EE\left[f(Y_n)| Y_n\­0\right]=\theta_n^{-1} \EE\left[f(Y_{n,1})-f(Y_{n,1}^{(2:r_n)}) | X_{n,1}\­0\right] + o(1),$$
where $o(1)$ converges to $0$ as $n\to \infty$ uniformly for all bounded cluster functionals $f\in\mc{F}$, and 
$$\theta_n:= \dfrac{\PP\{Y_n\­ 0\}}{r_n v_n}=\PP\{ Y_{n,1}^{(2: r_n)}=0| X_{n,1}\­ 0\}(1+o(1)).$$
Additionally, if the assumption (TC) is satisfied, then:
$$
\begin{array}{cll}
m_W &:=&\sup\{i\³1: W_i\­0\}< \infty,
\\ \theta_n& \underset{_{n\to \infty}}{\lfled}&\theta:=  \PP\{ W_i=0, \forall i\³ 2\}=\PP\{ m_W=1\}>0,
\\ P_{f(Y_n)| Y_n\­ 0}& \overset{w}{\underset{_{n\to \infty}}{\lfled}}& \dfrac{1}{\theta} \left(  \PP\{f(W)\in \cdot \}- \PP\{f(W^{(2:\infty)})\in \cdot , m_W\³ 2\}\right).
\end{array}
$$
\end{prop}
\begin{coro} Suppose that the hypotheses from Proposition 1 are maintained and that the assumption (TC) is satisfied. If, additionally, for each case of weak dependence $\theta$, $\eta$, $\kappa$ and $\lambda$ we request that $r_n^2=\mc{O}(n^2 v_n^{p+1})$, $r_n^4=\mc{O}(n^3 v_n^{p+1})$, $r_n^4=\mc{O}(n^3 v_n^{2p+1})$ and $r_n^4=\mc{O}(n^3 v_n^{2p+1})$ be fulfilled respectively for some $p>0$, then the fidis of the cluster functionals empirical process $(Z_n(f))_{f\in\mc{F}}$ converge to the fidis of a centered Gaussian process $(Z(f))_{f\in\mc{F}}$ with covariance function $c$ defined by
\begin{align}\label{cc}
c(f,g)=\EE\left[(fg)(W)- (fg)(W^{(2:\infty)})\right].
\end{align}
\end{coro}
\section{Application: the extremogram}
Hereafter we shall provide an application where it is enough to consider the  convergence finite - dimensional of the EPCFs. Specifically, we will prove that under suitable distributional conditions, the \citet{Davis2009}'s extremogram estimator for weakly dependent time series is asymptotically normal.

\subsection{The extremogram}
For a strictly stationary $\RR^d$ - valued time series $(X_t)_{t\in\ZZ}$, \citet{Davis2009} have defined the \textbf{extremogram} the two sets $A$ and $B$ bounded away from zero\footnote{A set $S$ is bounded away from zero if $S\subset \{y: |y|>r\}$ for some $r>0$} by 
\begin{align}\label{extrem}
\rho_{A,B}(h):=\lim_{x\to\infty} \PP\{ x^{-1} X_h\in B| x^{-1} X_0 \in A\} \qquad h=0,1, 2,  \ldots
\end{align}
provided the limit exist. 

As \citet{Davis2009} have said, a "natural" estimator of the extremogram based on the observations $X_1, \ldots, X_n$ is:
\begin{align}\label{estimator}
\h{\rho}_{A,B,n}(h):=\dfrac{\sum_{i=1}^{n-h} \mb{1} \{ u_n^{-1}X_{i+h}\in B, u_n^{-1} X_i \in A\}}{\sum_{i=1}^n \mb{1} \{u_n^{-1} X_i \in A\}},
\end{align}
where $u_n$ is a high quantile of the process which replaces $x$ in the limit (\ref{extrem}). Of course, the choice of such a sequence of quantiles $(u_n)_{n\in\NN}$ is not arbitrary. Particularly, such a sequence must satisfy the following condition\footnote{A sufficient condition permitting the Condition (\ref{RV})  to be fulfilled, involves the process to be regularly varying with index $\alpha>0$. For more details on the interpretation of the structure of regularly varying sequences, see \citep{Basrak2009}.}: 
\begin{align}\label{RV}
n \PP\{ u_n^{-1}(X_1, \ldots, X_h)\in \cdot )  \underset{n\to\infty}{\overset{vague}{\lfled }} \mu_h(\cdot), \end{align}
for each $h\³ 1$, where $(\mu_h)_{h\in\NN}$ is a sequence of non-null Radon measures on the Borel $\sigma$-field of $\RR^{dh}\setminus \{0\}$. \\Besides,  $v_n=\PP\{ u_n^{-1}X_0  \in  A\} \underset{n\fled\infty}{\lfled} 0$ with $n v_n\underset{n\fled\infty}{\lfled} \infty$ in order to have consistency in the results.
 
Let us now define the pre-asymptotic extremogram (PA-extremogram) $\rho_{A,B,n}(h):=\PP\{ u_n^{-1} X_h\in B| u_n^{-1} X_0 \in A\}$ and let $l$ be a positive integer. Then, under suitable conditions of convergence and weak dependence\footnote{Under suitable $\alpha$-mixing conditions, \citet{Davis2009, Davis2012} proved the convergence (\ref{extconv}).},
\begin{align}\label{extconv}
\sqrt{n v_n}\left( \h{\rho}_{A,B,n}(h) - \rho_{A,B,n}(h)\right)_{0\²h\²l}  \underset{n\to\infty}{\overset{\mc{D}}{\lfled }}   \mc{N}\left(0, \Sigma_{A,B}\right), 
\end{align} 
where $\Sigma_{A,B}$ is defined in (\ref{matrixcova}) below. \\Indeed, if for each $h\in\{1,\ldots, l\}$ with $l<r$, we define the cluster functional  $f_{A,B,h} : (\RR^d_\cup, \mc{R}_\cup)\lfled (\RR, \mc{B}(\RR))$  such that
\begin{align}\label{funcext}
f_{A,B,h}(x_1,\ldots, x_r):=\sum_{i=1}^{r-h} \1 \{ x_i\in A, x_{i+h}\in B\},
\end{align}
then, by using the normalization (\ref{norextremogram}), we can rewrite the estimator (\ref{estimator})  as: 
\begin{align}\label{rhoestimator}
\h{\rho}_{A,B,n}(h)=\dfrac{\sqrt{nv_n} Z_n(f_{A,B,h}) + m_n \EE f_{A,B,h}(Y_{n,1}) + \sum_{j=1}^{m_n} \delta_{n,j} (f_{A,B,h}) + R_n(A,B,h)}{\sqrt{nv_n} Z_n(f_{A,A,0}) + m_n \EE f_{A,A,0}(Y_{n,1}) + \sum_{j=1}^{m_n} \delta_{n,j} (f_{A,A,0}) + R_n(A,A,0)}
\end{align}
where 
\begin{align}
\delta_{n,j} (f_{A,B,h}):&=\sum_{i=j r_n - h +1}^{j r_n} \1 \{ u_n^{-1} X_i \in A, u_n^{-1} X_{i+h} \in B\}
\\R_n(A, B, h):&=\sum_{i=m_n r_n +1 }^{n-h}\1\{ u_n^{-1} X_i \in A, u_n^{-1} X_{i+h} \in B\}.
\end{align}
Then, if $(X_{n,i})_{1\²i\² n, n\in\NN}$ satisfies any of the following weak dependency conditions: 
\begin{enumerate}
\item[(D.1'')] $\theta$ (resp. $\eta$)-weakly dependent such that $\sum_{k=1}^{n-1} \theta_n(k)$ (resp. $\sum_{k=1}^{n-1} \eta_n(k)$) $= \mc{O}(v_n^{p+1})$ 
\item[(D.2'')] $\kappa$ ( resp. $\lambda$) -weakly dependent such that $\sum_{k=1}^{n-1} \kappa_n(k)$   (resp. $\sum_{k=1}^{n-1} \lambda_n(k))$ $= \mc{O}(v_n^{2p+1})$,
\end{enumerate}
for some $p>0$, and
\begin{align}\label{Segundo}
\EE\left( \sum_{i=1}^{r_n} \1 \{ u_n^{-1}X_i \in A\}\right)^2=\mc{O}(r_n v_n)
\end{align}
where $r_n \ll v_n^{-1} \ll n$  and $\sqrt{n v_n}=o(r_n)$, with $r_n \underset{n\fled\infty}{\lfled}\infty$. Then, 
\begin{align}\label{asympt}
\sqrt{n v_n}\left(  \h{\rho}_{A,B,n}(h) - \rho_{A,B,n}(h)\right) = Z_n(f_{A,B,h})- \rho_{A,B,n}(h) Z_n(f_{A,A,0}) + o_P(1).
\end{align}
Now, based on the equality (\ref{asympt}), we shall formalize (\ref{extconv})  through the following result: 
\begin{prop}
Assume that the big and small blocks sizes $r_n$ and $l_n$ are such that $l_n \ll r_n \ll v_n^{-1} \ll n$  and $\sqrt{n v_n}=o(r_n)$, with $l_n \underset{n\fled\infty}{\lfled}\infty$. Moreover, suppose that the normalization (\ref{norextremogram}), built from a stricly stationary regularly varying sequence $(X_i)_{i\in\NN}$ of  $\RR^d$-valued random vectors, satisfies:
\begin{enumerate}
\item[(i)]  some of the weak dependence conditions of the list (D.1") - (D.2");
\item [(ii)] the concentration condition (\ref{Concent}); 
\item [(iii)] the convergence conditions: (C.3) with $f=f_{A,B,h} $, and 
 \begin{align}\label{feacondition}
 \left( \EE \left| f_{A,B,h}(Y_n) \right| ^{2+\delta}\right)^{1/2} =o\left( \dfrac{r_n^2 (n v_n)^{(\delta+2)/4}}{n^2}\right), 
 \end{align}
for some $\delta\in (0,6];$ and that  
 \item[(iv)] there exists the covariance functions: 
 \begin{align}
 (r_n v_n)^{-1}\sum_{i=1}^{r_n-h}\sum_{j=1}^{r_n-h'} &\PP\{ X_{n,i}, X_{n,j}\in A; X_{n,i+h}, X_{n, j+h'}\in B\} \underset{n\fled\infty}{\lfled}  \sigma_{A,B}(h,h')\label{cova1}
 \\
(r_n v_n)^{-1}\sum_{i=1}^{r_n-h}&\sum_{j=1}^{r_n} \PP\{ X_{n,i}, X_{n,j}\in A; X_{n,i+h}\in B\} \underset{n\fled\infty}{\lfled}  \sigma'_{A,B}(h)\label{cova2}.
 \end{align}
 \end{enumerate}
Then, 
$$\sqrt{n v_n}\left( \h{\rho}_{A,B,n}(h) - \rho_{A,B,n}(h)\right)_{0\²h\²l}  \underset{n\to\infty}{\overset{\mc{D}}{\lfled }}   \mc{N}\left(0, \Sigma_{A,B}\right),$$
where the covariance matrix is defined by
\begin{align}\label{matrixcova}
\Sigma_{A,B}:=\left[ \sigma_{A,B}(h,h')-\rho_{A,B}(h') \sigma'_{A,B}(h) - \rho_{A,B}(h) \sigma'_{A,B}(h')+\rho_{A,B}(h)\rho_{A,B}(h')\sigma'_{A,A}(0)\right]_{0\² h,h'\² l}.
\end{align}
\end{prop}
\subsection{Simulation Study}
In order to numerically determine our results, we shall do a small simulation of a real-valued weakly dependent data, estimate its extremogram by means of (\ref{estimator})  and compare it to the theoretic pre-asymptotic extremogram $\rho_{A,B,n}(\cdot)$. 
\subsubsection{Theoretic model} Let us consider the AR(1)-process (\ref{AR}) given in the introduction. Here, as $X_0$ is uniformly distributed on $[0,1]$ and $X_i=b^{-i} X_0 + \sum_{s=1}^i b^{s-i-1} \xi_s$ for all $i\³1$, then for $A=B=(1,\infty)$ we obtain:
\begin{align}\label{ExtremogramAR1}
\PP\left\{ \dfrac{X_h}{y}\in B\left| \dfrac{X_0}{y}\in A\right.\right\}=\dfrac{1}{b^h} \sum_{j_1, \ldots, j_h\in U(b)}\min \left\{ 1, \dfrac{1}{1-y} \left(1-yb^h+\sum_{s=1}^h \dfrac{j_s}{b^{1-s}} \right)_+\right\},
\end{align}
where $y=1-1/x$. Therefore, 
\begin{align}\label{existe}
\rho_{A,B}(h) &= b^{-h} \qquad \text{for }  h=0,1,\ldots ; \text{ and }
\\
\rho_{A,B,n}(h) &= \dfrac{1}{b^h} \sum_{j_1, \ldots, j_h\in U(b)}\min \left\{ 1, \dfrac{1}{v_n} \left(1-(1-v_n)b^h+\sum_{s=1}^h \dfrac{j_s}{b^{1-s}} \right)_+\right\},\label{PA}
\end{align}
for $h< n$, where $v_n:=\PP\{ X_0/u_n \in A\}=1-u_n$.

Note that the equality (\ref{existe}) proves that the family of weakly dependent processes such that the extremogram (\ref{extrem}) exists is not empty.\footnote{\citet{Davis2009} have proved under suitable mixing conditions and distributional assumptions that the limit in (\ref{extrem}) exists.}.
\subsubsection{Experiment} In order to carry out the experiment, we have generated from the AR(1)-process defined in (\ref{AR}) with $b=2$, $N=50$ samples of $n=2000$. Besides, here we have particularly taken $v_n=1/\sqrt{n}=1/10\sqrt{2}$. Therefore, the high quantile is $u_n=(10\sqrt{2}-1)/10\sqrt{2}$, assuming that we have determined that the distribution of the $X_i$ is uniformly distributed on $[0,1]$.

The mean PA-extremogram and the theoretic PA-extremogram (\ref{PA}) for lags $h=1\ldots, 20$ corresponding to the right tail $(A=B=(1,\infty))$ are displayed as the blue solid line and the black solid line, respectively,  in the left panel of Fig.1.  In the same panel, we show $95\%$ confidence bands (red dashed lines), symmetric with respect to the mean PA-extremogram. As was expected for the AR(1)-process studied here, when observing the confidence bands, note that the extremal dependence vanishes as the lag $h$ increases.

In the right panel of Fig.1, we show $95\%$ confidence bands (red dashed lines) for the statistical errors $\e_h(x_1^j, \ldots, x_n^j):=\h{\rho}_{A,B,n}(h, (x_1^j, \ldots, x_n^j))-\rho_{A,B,n}(h)$, for $j=1,\ldots, N$, with $\rho_{A,B,n}(h)$ defined in (\ref{PA}). {\color{blue} }
\begin{figure}
\centering
\begin{minipage}{.5\columnwidth}
\centering
\vspace{-0.5cm}
\includegraphics[width=\columnwidth]{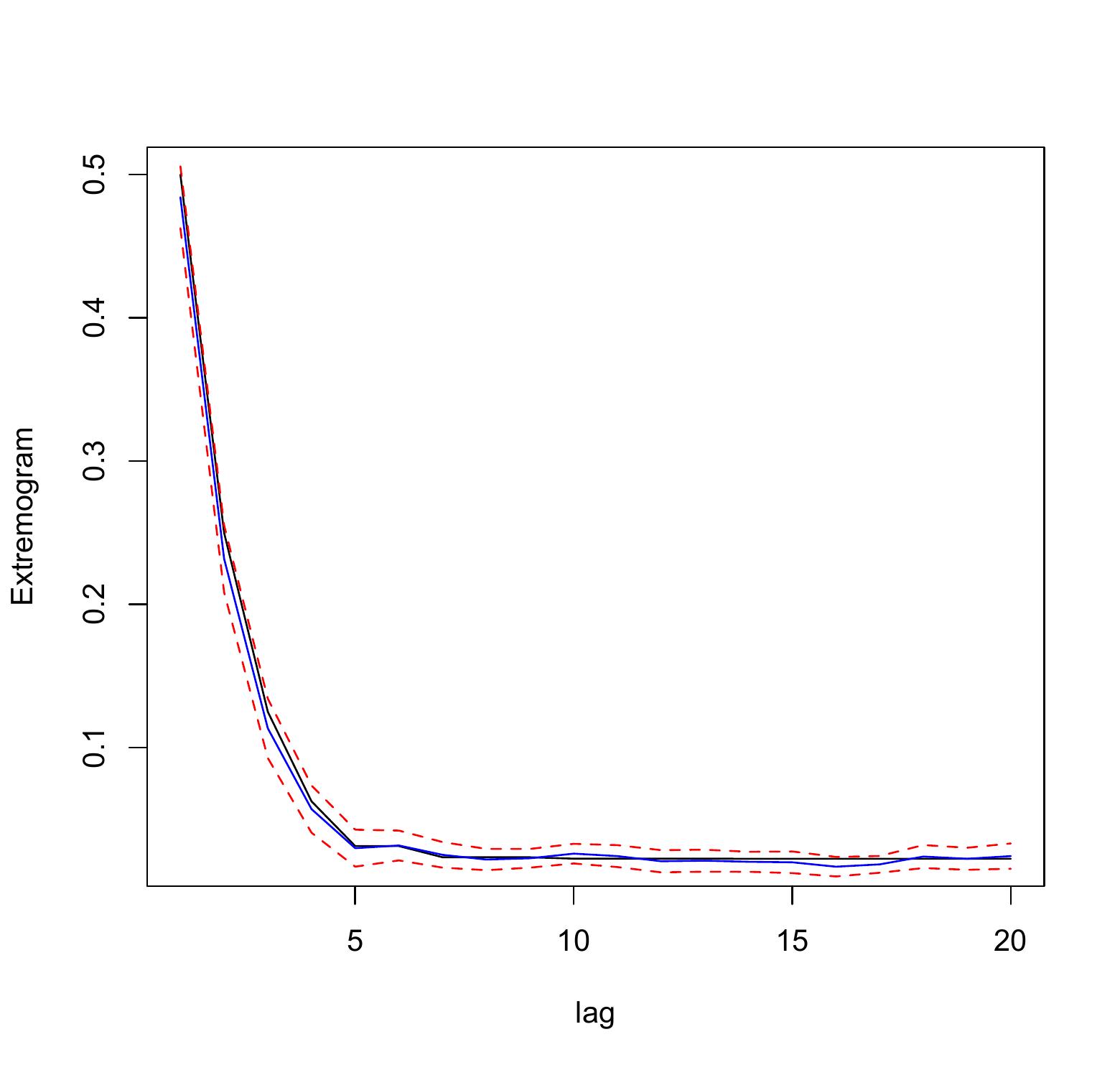}
\end{minipage}%
\begin{minipage}{.5\columnwidth}
\centering
\includegraphics[width=\columnwidth]{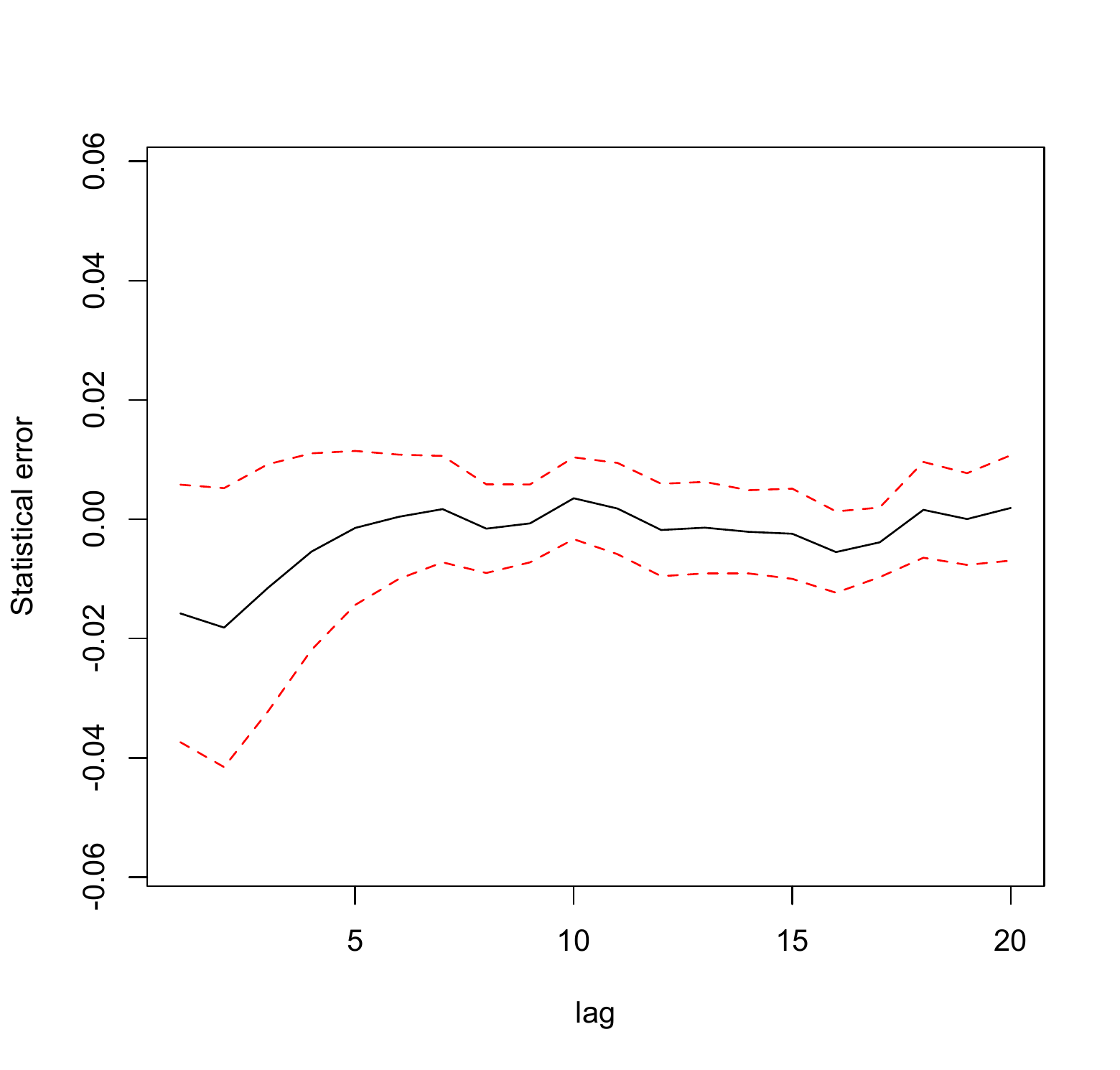}
\label{fig:b3}
\end{minipage}
\textbf{Fig.1.} \small Left: $95\%$ confidence bands (red dashed lines) for the PA-extremogram (black solid line) of the AR(1)-process (\ref{AR}) with $b=2$, and the mean estimated PA-extremogram  (blue solid line). Right: $95\%$ confidence bands (red dashed lines) for the statistical errors of the PA-extremogram and the mean statistical error of the  PA-extremogram (black solid line).
\end{figure}
\section{Proofs}\
{\it \textbf{Proof of Theorem 1.}}
The proof is basically a direct application of  Theorem 1 - \citep{Bardet2007} to the random variables $(W_{n,j})$ defined in (\ref{N2}).

First, notice that Assumption (C.1) implies that $$B_n(\epsilon):=\sum_{j=1}^{m_n} \mb{E}\|W_{n,j}\|^2\mb{1}\left\{\|W_{n,j}\|>\epsilon \right\}\underset{n\fled\infty}{\lfled}0,$$ for all k-tuple $(f_1, \ldots, f_k)\in\mc{F}^k$, $k\in\NN$ and $\epsilon>0$. Note that this last statement is weaker than Assumption $H_\delta$ of \citet{Bardet2007}, required in the assumptions of their theorem. \\On the other hand, Assumption (C.2) ensures the existence of the positive matrix $\Sigma_k=(c(f_i,f_j))_{i,j=1,\ldots,k}$ such that $$\Sigma_{n,k}:=(\Cov(f_i(Y_n),f_j(Y_n)))_{i,j=1,\ldots,k} \underset{n\fled\infty}{\lfled} \Sigma_k,$$ for all k-tuple $(f_1, \ldots, f_k)\in\mc{F}^k$, $k\in\NN$. The proof ends considering the condition of dependence (\ref{TDep}). 
\hfill$\square$\vskip2mm\hfill
\\
{\it \textbf{Proof of Proposition 1.}} 
The proof of this proposition is based on Theorem 1. Therefore, we only have to prove that $T_t(m_n| f_1,\ldots, f_k) \underset{n\fled\infty}{\lfled} 0$, for all $t\in\RR^k$ and for all $k$-tuple of cluster functionals $(f_1,\ldots, f_k)\in\mc{F}^k$. Indeed, for $j\in\{2,\ldots, m_n\}$ and  $(f_1,\ldots, f_k)\in\mc{F}^k$, notice that
$$\Cov\left(\e^{i<t, \sum_{s=1}^{j-1}W_{n,s}(f_1, \ldots, f_k)>}, \e^{i<t,W_{n,j}(f_1, \ldots, f_k)>}\right)$$
can be rewritten as: 
\begin{align}\label{l1}
\Cov(F_j, G_j):=\Cov\left(F_{t,n}^{(f_1,\ldots, f_k)}(Y_{n,1},\ldots, Y_{n,j-1}), G_{t,n}^{(f_1, \ldots, f_k)}(Y_{n,j})\right),
\end{align}
where $G_{t,n}^{(f_1, \ldots, f_k)}(s)=\e^{i<t,\sum_{l=1}^k\frac{f_l(s)-\EE f_l (s)}{\sqrt{n v_n}}e_l>}$; $F_{t,n}^{(f_1,\ldots, f_k)}(s_1, \ldots, s_j)=\prod_{h=1}^j G_{t,n}^{(f_1,\ldots, f_k)}(s_h)$ and ${e_1, \ldots, e_k}$ is the canonical base in $\RR^k$. Moreover, it is clear that $\| G_{t,n}^{(f_1, \ldots, f_k)}\|_\infty\²1$ and $\| F_{t,n}^{(f_1, \ldots, f_k)}\|_\infty\²1$, for all $t\in \RR^k$, $\vec{f}=(f_1, \ldots, f_k)\in\mc{F}^k$, for any $k\³1$. Then,
\begin{align}
\left|\Cov\left(F_j,G_j\right)\right|&=\left|\Cov\left(F_{t,n}^{\vec{f}}(Y_{n,1},\ldots, Y_{n,j-1}), G_{t,n}^{\vec{f}}(Y_{n,j})\right)\right|
\nonumber\\&\²\left|\Cov\left(F_{t,n}^{\vec{f}}(Y_{n,1},\ldots, Y_{n,j-1})-F_{t,n}^{\vec{f}_\rho^{[T]}}(Y_{n,1}^{(r_n-l_n)},\ldots, Y_{n,j-1}^{(r_n-l_n)}), G_{t,n}^{\vec{f}}(Y_{n,j})\right)\right|
\nonumber\\&+\left|\Cov\left(F_{t,n}^{\vec{f}_\rho^{[T]}}(Y_{n,1}^{(r_n-l_n)},\ldots, Y_{n,j-1}^{(r_n-l_n)}), G_{t,n}^{\vec{f}}(Y_{n,j})-G_{t,n}^{\vec{f}_\rho^{[T]}}(Y_{n,j}^{(r_n-l_n)})\right)\right|
\nonumber\\&+\left|\Cov\left(F_{t,n}^{\vec{f}_\rho^{[T]}}(Y_{n,1}^{(r_n-l_n)},\ldots, Y_{n,j-1}^{(r_n-l_n)}), G_{t,n}^{\vec{f}_\rho^{[T]}}(Y_{n,j}^{(r_n-l_n)})\right)\right|
\nonumber\\ &\²2\mb{E}\left|F_{t,n}^{\vec{f}}(Y_{n,1},\ldots, Y_{n,j-1})-F_{t,n}^{\vec{f}_\rho^{[T]}}(Y_{n,1}^{(r_n-l_n)},\ldots, Y_{n,j-1}^{(r_n-l_n)})\right|
\label{l3}\\&+2\mb{E}\left|G_{t,n}^{\vec{f}}(Y_{n,j})-G_{t,n}^{\vec{f}_\rho^{[T]}}(Y_{n,j}^{(r_n-l_n)})\right|
\label{l4}\\&+\left|\Cov\left(F_{t,n}^{\vec{f}_\rho^{[T]}}(Y_{n,1}^{(r_n-l_n)},\ldots, Y_{n,j-1}^{(r_n-l_n)}), G_{t,n}^{\vec{f}_\rho^{[T]}}(Y_{n,j}^{(r_n-l_n)})\right)\right|
\label{l5},
\end{align}
where $\vec{f}^{[T]}_\rho=(f_{1,\rho}^{[T]},\ldots, f_{k,\rho}^{[T]})$ is a $k$-tuple of Lipschitzian cluster functionals which approximates $\vec{f}$ as $\rho\downarrow 0$ and as $T\uparrow \infty$, defined as follows.

Let $f: E^r\lfled \RR$ be a cluster functional. First, we consider $f^{[T]}=f\vee(-T)\wedge T$, a truncation of $f$ by $T$, for $T>0$. Now, we define the set
\begin{align}
D^r(f^{[T]})=\{y\in E^r: \text{ either } y \text{ is a discontinuity point of } f^{[T]} \text{ or } \|D^+ f^{[T]}(y)\|=\infty\}
\nonumber,
\end{align}
where $D^+ f(y)$ denotes the upper Dini derivative matrix at $y$ of $f:E^r \lfled \RR$. 
So, for each $\rho>0$ we denote
$$D_\rho^r=\bigcup \left\{B(y, \rho/2): y\in D^r(f^{[T]})\right\}.$$Clearly $D_\rho^r$ is open in $E^r$, thus $C_\rho^r=E^r\setminus D_\rho^r$ is closed in $E^r$. Moreover, $L=\RR$ is an affine space of ``type $m$" \footnote{An affine space if of type $m$ if for each first countable space $X$ and any every continuous $f: X\lfled L$, we have that for each $x\in X$ and nbd $W\supset f(X)$, there exists a nbd $U$ of $x$ and a convex set $C\subset L$ such that $f(U)\subset C \subset W$.\\ Thus, all locally convex linear topological spaces (and also all vector spaces with the finite topology) are of type $m$. For more details see \citep{Dugundji}.}. Therefore, using Theorem 6.1.-  \citep{Dugundji}, the continuos functional $\left.f^{[T]}\right|_{C_\rho}: C_\rho\lfled \RR$ has a continuous extension $g_\rho:E^r \fled \RR$ such that $g_\rho(E^r)\subset [\text{convex hull of } f^{[T]}(C_\rho)]$. Finally, we can choose $f_\rho^{[T]}=g_\rho$ and note that
\begin{align}\label{l6}
\Lip f^{[T]}_\rho\² \dfrac{2T}{\rho r\sqrt{d}}.
\end{align}
Therefore, we can easily obtain bounds for $G^{\vec{f}}_{t,n}$ and $F^{\vec{f}}_{t,n}$:
\begin{align}\label{Lip1}
\Lip G^{\vec{f}}_{t,n}\² \dfrac{2 \|t\|}{\sqrt{n v_n}}\sqrt{\sum_{l=1}^k \left(\Lip f_{l,\rho}^{[T]}\right)^2 }\² \dfrac{4 \|t\| \sqrt{k} T}{\rho r\sqrt{d}\sqrt{n v_n}} \qquad \text{ and } \qquad \Lip F^{\vec{f}}_{t,n}\² \Lip G^{\vec{f}}_{t,n}.
\end{align}
Denote by 
$$C(F,G):=\sum_{j=1}^{m_n}\left|\Cov\left(F_{t,n}^{\vec{f}^{[T]}_\rho}(Y_{n,1}^{(r_n-l_n)},\ldots,Y_{n,j-1}^{(r_n-l_n)}),G_{t,n}^{\vec{f}^{[T]}_\rho}(Y_{n,j}^{(r_n-l_n)})\right)\right|,$$
the sum over $j=1,\ldots, m_n$ of the term (\ref{l5}). Now, combining (\ref{Lip1}) with the definition of the weak dependence coefficients (\ref{dp1}) - (\ref{dp4}), we obtain bounds for $C(F,G)$ according to the respective condition of weak dependence assumed for $(X_{n,i})_{1\²i\²n, n\in\NN}$:
\begin{enumerate}
\item $\theta$ - WD implies $C(F,G)\²\dfrac{4T\|t\|\sqrt{k}}{\rho \sqrt{d}}\dfrac{n}{r_n \sqrt{n v_n}}\theta_n(l_n)$
\item $\eta$ - WD implies $C(F,G)\²\dfrac{2 T \|t\| \sqrt{k}}{\rho \sqrt{d}}\dfrac{n^2}{r_n^2 \sqrt{n v_n}}\left(1+\dfrac{r_n}{n}\right)\eta_n(l_n)$
\item $\kappa$ - WD implies $C(F,G)\²\dfrac{8 T^2 \|t\|^2 k}{\rho^2 d}\dfrac{n}{r_n^2v_n} \kappa_n(l_n)$
\item $\lambda$ - WD implies $C(F,G)\²\left[ \dfrac{2T\|t\|\sqrt{k}}{\rho \sqrt{d}}\dfrac{n^2}{r_n^2 \sqrt{n v_n}}\left(1+\dfrac{r_n}{n}\right) + \dfrac{8T^2 \|t\|^2k}{\rho^2 d}\dfrac{n}{r_n^2 v_n}\right]\lambda_n(l_n)$,
\end{enumerate}
For the sum of the terms (\ref{l3}) and (\ref{l4}), notice that
\begin{align}\label{l7}
|F_{t,n}^f(s_1, \ldots, s_p)-F_{t,n}^{f'}(s'_1, \ldots, s'_p)|\²& \sum_{i=1}^p|G_{t,n}^f(s_i)-G_{t,n}^f(s'_i)|+\sum_{i=1}^p|G_{t,n}^f(s'_i)-G_{t,n}^{f'}(s'_i)|.
\nonumber
\end{align} 
Thus, we develop the sum of the terms (\ref{l3}) and (\ref{l4}):
\begin{align}
&2\mb{E}\left|F_{t,n}^{\vec{f}}(Y_{n,1},\ldots, Y_{n,j-1})-F_{t,n}^{\vec{f}_\rho^{[T]}}(Y_{n,1}^{(r_n-l_n)},\ldots, Y_{n,j-1}^{(r_n-l_n)})\right|+2\mb{E}\left|G_{t,n}^{\vec{f}}(Y_{n,j})-G_{t,n}^{\vec{f}_\rho^{[T]}}(Y_{n,j}^{(r_n-l_n)})\right|
\nonumber\\ \²&2\EE\left(\sum_{i=1}^{j-1}|G_{t,n}^{\vec{f}}(Y_{n,i})-G_{t,n}^{\vec{f}}(Y_{n,i}^{(r_n-l_n)})|+\sum_{i=1}^{j-1}|G_{t,n}^{\vec{f}}(Y_{n,i}^{(r_n-l_n)})-G_{t,n}^{\vec{f}_\rho^{[T]}}(Y_{n,i}^{(r_n-l_n)})|\right)
\nonumber\\+&2\EE\left(|G_{t,n}^{\vec{f}}(Y_{n,j})-G_{t,n}^{\vec{f}}(Y_{n,j}^{(r_n-l_n)})|+|G_{t,n}^{\vec{f}}(Y_{n,i}^{(r_n-l_n)})-G_{t,n}^{\vec{f}_\rho^{[T]}}(Y_{n,i}^{(r_n-l_n)})|\right)
\nonumber
\end{align}
\begin{align}
=&2\sum_{i=1}^{j}\EE|G_{t,n}^{\vec{f}}(Y_{n,i})-G_{t,n}^{\vec{f}}(Y_{n,i}^{(r_n-l_n)})|+2 \sum_{i=1}^{j}\EE|G_{t,n}^{\vec{f}}(Y_{n,i}^{(r_n-l_n)})-G_{t,n}^{\vec{f}_\rho^{[T]}}(Y_{n,i}^{(r_n-l_n)})|
\nonumber\\ =& 2j\EE|G_{t,n}^{\vec{f}}(Y_n)-G_{t,n}^{\vec{f}}(Y_n^{(r_n-l_n)})|+2j\EE|G_{t,n}^{\vec{f}}(Y_n^{(r_n-l_n)})-G_{t,n}^{\vec{f}_\rho^{[T]}}(Y_n^{(r_n-l_n)})|
\nonumber\\ =& 2j \|t\| \sqrt{\sum_{l=1}^k \dfrac{\Var(\Delta_n(f_l))}{n v_n}}+ 2j\|t\| \sqrt{\sum_{l=1}^k \dfrac{\Var(f_l(Y_n^{(r_n-l_n)})-f_{l,\rho}^{[T]}(Y_n^{(r_n-l_n)}))}{n v_n}}
\nonumber\\ \²&\dfrac{2j\|t\|}{\sqrt{n v_n}}\left(\sqrt{\sum_{l=1}^k\Var(\Delta_n(f_l))}+\sqrt{\sum_{l=1}^k\EE|f_l(Y_n^{(r_n-l_n)})|^2\mb{1}\{|f_l(Y_n^{(r_n-l_n)})|>T\}+ 4T^2k (C \rho^\alpha)^{(r_n-l_n)\delta}}\right)
\nonumber,
\end{align}
for some $\delta\in(0,1)$.

This proves 
\begin{align}
T_t(m_n|f_1,\ldots, f_k)&=\sum_{j=1}^{m_n} |\Cov(F_j,G_j)|
\nonumber\\ &\² \|t\| \sqrt{k}\left(1+\dfrac{r_n}{n}\right)\dfrac{n^2}{r_n^2 \sqrt{n v_n}}\left( \sqrt{\Var(\Delta_n(f))}\right.
\\&+ \sqrt{ \EE|f(Y_n^{(r_n-l_n)})|^2\mb{1}\{|f(Y_n^{(r_n-l_n)})|>T\}}
\\&+\left. 2 T (C \rho^\alpha)^{(r_n-l_n)\delta/2}\right)+C(F,G).
\end{align}

Finally it suffices to choose $T$ such that $T=\mc{O}(\sqrt{n v_n})$ with 
$$ \text {1) } \rho=\left(\dfrac{r_n}{n}\right)^{\frac{2}{2+\alpha r_n \delta}}\dfrac{(\theta_n(l_n))^{\frac{2}{2+\alpha \delta r_n}}}{C^{\frac{1}{\alpha + \frac{2}{\delta r_n}}}}, \text{  2) } \rho=\left(\dfrac{r_n}{n}\right)^{\frac{2}{2+\alpha r_n \delta}}\dfrac{(\eta_n(l_n))^{\frac{2}{2+\alpha \delta r_n}}}{C^{\frac{1}{\alpha + \frac{2}{\delta r_n}}}},$$ 
$$ \text{ 3) } \rho=\left(\dfrac{r_n}{n}\right)^{\frac{2}{4+\alpha r_n \delta}}\dfrac{(\kappa_n(l_n))^{\frac{2}{4+\alpha \delta r_n}}}{C^{\frac{1}{\alpha + \frac{4}{\delta r_n}}}}, \text{ 4) } \rho=\left(\dfrac{r_n}{n}\right)^{\frac{2}{4+\alpha r_n \delta}}\dfrac{(\lambda_n(l_n))^{\frac{2}{4+\alpha \delta r_n}}}{C^{\frac{1}{\alpha + \frac{4}{\delta r_n}}}},$$ for each respective weak - dependence condition and take into account the assumptions (C.3) and (C.4) to obtain the CLT.
\hfill$\square$\vskip2mm\hfill
\\
{\it \textbf{Proof of Proposition 2.}} 
Suffices to prove the following multidimensional version of \citet{Segers2003}'s condition (6) 
\begin{align}\label{CondiDepend}
\lim_{l\to\infty}\limsup_{n\to\infty} \PP \left\{\left. Y_n^{(l+1:r_n)} \­ 0\right| X_{n,1}\­0 \right\}=0,
\end{align}
since the rest of the proof follows the same steps of the proof of \citet{Drees2010}'s Lemma 2.5.
 
Indeed, let $f(\cdot)=\1\{\cdot \­ 0\}$ and $g(\cdot)=\1\{\cdot \­ 0\}$ be functions defined on  $E^{r_n-l}$ and $E$ respectively. Now, consider increasing sequences of functions $f_k: E^{r_n-l}\lfled [0,1]$ and $g_k:E\lfled [0,1]$ which approximate to $f$ and $g$ respectively; and such that $\Lip (f_k)=\Lip(g_k)=v_k^{-p}$ for some $p>0$. 
 Of course, note that $k=k(n)  \ll n$. \\ 
Then, if the random variables $(X_{n,i})_{1\²i\²n}$ are $(\epsilon_n, \psi)$-WD we have that
\begin{multline}
\limsup_{n\to\infty}\PP \left\{\left. Y_n^{(l+1:r_n)} \­ 0\right| X_{n,1}\­0 \right\}=\limsup_{n\to\infty} v_n^{-1} \EE\1\{ Y_n^{(l+1:r_n)}\­0\}\1\{X_{n,1}\­0\}
\\=\limsup_{k\to\infty} v_n^{-1}  \EE f_k(Y_n^{(l+1:r_n)}) g_k(X_{n,1})
\\= \limsup_{k\to\infty} v_n^{-1}  \Cov\left(f_k(Y_n^{(l+1:r_n)}), g_k(X_{n,1})\right)+ \limsup_{k\to\infty} v_n^{-1}    \EE f_k(Y_n^{(l+1:r_n)}) \EE g_k(X_{n,1})
\\ \² \limsup_{k\to\infty} \psi(r_n-l, 1, \Lip(f_k), \Lip(g_k)) \dfrac{\epsilon_n(l)}{v_n} + \limsup_{k\to\infty}  r_n v_n(1-\dfrac{l}{r_n})
\\ = \limsup_{k\to\infty}  \psi(r_n-l, 1, \Lip(f_k), \Lip(g_k)) \dfrac{\epsilon_n(l)}{v_n}
\end{multline}
Finally, if the random variables $(X_{n,i})_{1\²i\²n, n\in\NN}$ are weakly dependent in some cases of the list (D.1') - (D.3'), then the limit (\ref{CondiDepend}) is proven.
\hfill$\square$\vskip2mm\hfill
\\
{\it \textbf{Proof of Corollary 1.}} 
Note that if for each case of weak dependence  $\theta$, $\eta$, $\kappa$ and $\lambda$ we ask that $r_n^2=\mc{O}(n^2 v_n^{p+1})$, $r_n^4=\mc{O}(n^3 v_n^{p+1})$, $r_n^4=\mc{O}(n^3 v_n^{2p+1})$ and $r_n^4=\mc{O}(n^3 v_n^{2p+1})$ are fulfilled respectively, and we combine this with the conditions (D.1) and (D.2) of Proposition 1, then we have (D.1'), (D.2') and (D.3') respectively, as the weak dependence case may be.

\hfill$\square$\vskip2mm\hfill
\\
{\it \textbf{Proof of the relation (\ref{asympt}).}} 
Relation (\ref{rhoestimator}) implies that for each $h=0, \ldots, l$, 
\begin{multline*}
\sqrt{n v_n} (\h{\rho}_{A,B,n}(h)-\rho_{A,B,n}(h))
\\ = \dfrac{Z_n(f_{A,B,h}) - \frac{h \sqrt{n v_n}}{r_n}\rho_{A,B,n}(h) - \rho_{A,B,n}(h)Z_n(f_{A,A,0})+ S_n(h)+ D_n(h)}{(nv_n)^{-1/2}Z_n(f_{A,A,0}) + 1 + (n v_n)^{-1} \sum_{j=1}^{m_n} \delta_{n,j}(f_{A,A,0})+(nv_n)^{-1} R_n(A,A,0)} + o(1),
\end{multline*}
where $S_n(h):=(n v_n)^{-1/2} \left( R_n(A,B,h)-\rho_{A,B,n}(h) \cdot R_n(A,A,0)\right)$ and
$$D_n(h) := (n v_n)^{-1/2} \sum_{j=1}^{m_n} \left( \delta_{n,j}(f_{A,B,h})- \rho_{A,B,n}(h) \cdot \delta_{n,j}(f_{A,A,0}) \right).$$
Using Chebyshev's inequality on the random variables $(n v_n)^{-1/2} \sum_{j=1}^{m_n} \delta_{n,j}(f_{A,B,h})$ and $(n v_n)^{-1/2} R_n(A,B,h)$ followed by (\ref{Segundo}), we prove that these variables (and consequently, also  $S_n(h)$ and $D_n(h)$) converge to zero in probability. On the other hand, by using again Chebyshev's inequality on the random variable $\zeta_n=(nv_n)^{-1/2}Z_n(f_{A,A,0})$ combined with the stationarity of the time series, the approximation of the indicatrix function $f(\cdot)=\1 \{ \cdot \in A\}$ through lipschitzian increasing functions $f_k$ (as we did in the proof for Proposition 2) and the dependence condition of list (D.1'')-(D.2'') as the case may be, we can obtain the convergence to zero in probability of such a variable $\zeta_n$. Finally, given that $\sqrt{n v_n}=o(r_n)$, we obtain the relation (\ref{asympt}).
\hfill$\square$\vskip2mm\hfill
\\
{\it \textbf{Proof of Proposition 3.}} First, note that (\ref{feacondition}) implies (\ref{Segundo}). Therefore, since the restrictions on the size of the big and small blocks $r_n$ and $l_n$ are the same as in the previous proof and (D.1'')-(D.2'') are maintained, we obtain relation (\ref{asympt}). 

On the other hand, observe that (\ref{feacondition})  implies (C.4) for some $\delta>0$, and also (C.1) but for some $\delta\in (0,6]$. The existence of the covariance function $c$ of (C.2) for the functionals $f_{A,B,h}$ and $f_{A,A,0}$ is assumed through relation (\ref{cova1}). Besides, the concentration condition and the condition (C.3) are maintained just as in Proposition 2. 

Observe that the weak dependence conditions (D.1'')-(D.2'') are stronger than the dependence conditions (D.1)-(D.3). Therefore, considering the existence of covariance function (\ref{cova2}), we obtain the result. 
\hfill$\square$\vskip2mm\hfill
\\
\\
{\it \textbf{Proof of the expressions (\ref{ExtremogramAR1}), (\ref{existe}) and (\ref{PA}).}} 
Due to that $X_i=b^{-1} X_0 + \sum_{s=1}^i b^{s-i-1}\xi_s$ for all $i\³1$, then for $h\³0$ we have 
\begin{multline*}
\PP\left\{\dfrac{X_h}{y}\in B \left| \dfrac{X_0}{y}\in A\right.\right\}= \dfrac{1}{1-y} \PP\{ X_h> y, X_0>y\}
\\=\dfrac{1}{1-y} \PP\left\{ X_0> \max\left\{y,  \left( y- \sum_{s=1}^h \dfrac{\xi_s}{b^{1-s+h}}\right)b^h\right\}\right\}
\\= \dfrac{1}{(1-y) b^h} \sum_{j_1, j_2, \ldots, j_h \in U(b)}\PP\left\{ X_0> \max\left\{y,  \left( y- \sum_{s=1}^h \dfrac{j_s}{b^{1-s+h}}\right)b^h\right\}\right\}
\\= \dfrac{1}{b^h} \sum_{j_1, j_2, \ldots, j_h \in U(b)} \min\left\{1,  \dfrac{1}{1-y}\left(1-yb^h+\sum_{s=1}^h \dfrac{j_s}{b^{1-s}}\right)_+\right\}. 
\end{multline*}
This proves relation (\ref{ExtremogramAR1}). On the other hand, note that $\mu_b(j_1,\ldots, j_h):=1-b^h+\sum_{s=1}^h b^{s-1}j_s\² -1$ for all $(j_1,\ldots, j_h)\in U^h(b)\setminus \{(b-1, \ldots, b-1)\}$ and $\mu_b(b-1, b-1, \ldots, b-1)=0$. Therefore, by substituting  $y=1-1/x$ in (\ref{ExtremogramAR1}), and taking the limit when $x\to\infty$, we obtain (\ref{existe}). Finally, to prove  (\ref{PA}) it suffices to substitute $y=1-v_n$ in (\ref{ExtremogramAR1}). 
\hfill$\square$\vskip2mm\hfill



\end{document}